\documentclass{article}
\usepackage[numbers,sort&compress]{natbib}
\usepackage[table]{xcolor}
\usepackage{graphicx}
\usepackage{amsmath,amsfonts,amssymb}
\usepackage{graphicx}
\usepackage{multirow}
\usepackage{tikz,pgf}
\usepackage{url}
\usepackage{amsmath}
\usepackage{graphicx}
\usepackage{epsfig}
\usepackage{epstopdf}
\usepackage{color}
\usepackage{enumitem}
\usepackage{mathtools}
\def\disp{\displaystyle}

\def\tto{\;{\lower 1pt \hbox{$\rightarrow$}}\kern -10pt
	\hbox{\raise 2pt \hbox{$\rightarrow$}}\;}
\def\Hat{\widehat}
\def\Tilde{\widetilde}
\def\Bar{\overline}
\def\ra{\rangle}
\def\la{\langle}

\def\epsilon{\varepsilon}

\def\h{\hfill\Box}
\def\R{\mathbb{R}}

\def\oz{\bar{z}}

\def\s{\square}

\def\co{\text{\normalfont core}}
\def\core{\mbox{\rm core}\,}

\def\gph{\text{\normalfont gph}\,}

\def\epi{\mbox{\rm epi}\,}

\def\dom{\mbox{\rm dom}\,}

\def\rge{\mbox{\rm rge}\,}

\def\h{\hfill\square}

\def\emp{\emptyset}

\def\oR{\Bar{\R}}

\def\emp{\emptyset}

\def\oR{\Bar{\R}}

\setlist[enumerate,1]{itemsep=0.0ex,parsep=0.5ex,label={\rm(\alph*)},leftmargin=*, align=left}
\newcounter{lk}

\DeclarePairedDelimiter{\abs}{\lvert}{\rvert}

\usepackage[hidelinks]{hyperref}
\newtheorem{theorem}{Theorem}[section]
\newtheorem{proposition}[theorem]{Proposition}

\newtheorem{lemma}[theorem]{Lemma}
\newtheorem{corollary}[theorem]{Corollary}

\title{Duality Theory on Vector Spaces}

\author{Dang Van Cuong \thanks{ (corresponding author), Faculty of Computer Sciences, Vietnam - Korea University of Information and Communications Technology, The University of Da Nang, Da Nang, Vietnam (email: dvcuong@vku.udn.vn)} \and Tuyen Tran \thanks{Department of Mathematics and Statistics, Loyola University Chicago, Chicago, IL, 60660, USA (email: ttran18@luc.edu)}}

\date{}

\begin{document}
	\maketitle
\begin{abstract}
	In this paper, we study the Fenchel-Rockafellar duality and the Lagrange duality in the general frame work of vector spaces without topological structures. We utilize the geometric approach, inspired from its successful application by B. S. Mordukhovich and his coauthors in variational and convex analysis (see \cite{CBN21,CBNC,CBNG22,CBNG,m-book,mn-book}). After revisiting coderivative calculus rules and providing the subdifferential maximum rule in vector spaces, we establish conjugate calculus rules under qualifying conditions through the algebraic interior of the function's domains. Then we develop sufficient conditions which guarantee the Fenchel-Rockafellar strong duality. Finally, after deriving some necessary and sufficient conditions for optimal solutions to convex minimization problems, under a Slater condition via the algebraic interior, we then obtain a sufficient condition for the Lagrange strong duality.
\end{abstract}
\textbf{Keywords:} Algebraic interior, Algebraic core, Fenchel duality, Lagrangian duality, Calculus rules, Optimality conditions \\ \\
\textbf{Mathematics Subject Classification (2000):} 49J52, 49J53, 90C31

	\section{Introduction}
	
	Duality theory is a fundamental aspect in the theory of optimization and its applications. Given an optimization problem (\textit{primal}), it can be paired with a corresponding \textit{dual} problem. In many situations it is easier to solve the dual than the corresponding primal problem, which motivates the theory. There are two main duality methods in optimization, the  \emph{Fenchel duality} and the \emph{Lagrangian  duality}, for which there has been extensive research in both theory, and applications \cite{Bertsekas2003,Bauschke2011,bl,bl2,Borwein2000,bwsig,blreport1,bot1,bot-revis,bot2,bot3,CBNG,CBNG22,durea,g,Gow-Te-90,AG,HU1,HU2,JW,mn-book,KFNg,zduality}.
	
	The first branch of duality theory in optimization studies \emph{Fenchel duality}. The main question of the Fenchel duality problem is ``{\bf When does the Fenchel strong duality take place?}'' which is answered by the Fenchel-Rockafellar theorem (in $\R^n$ see \cite[Theorem 31.1]{r}). Recently, under qualifications on generalized interiors, some sufficient conditions have been given to get Fenchel strong duality on locally convex topological vector spaces; see, e.g., \cite{bl,bl2,bwsig,bot1,bot-revis,bot2,CBNG,CBNG22,zduality}. The primary objective of this paper is generalizing the Fenchel–Rockafellar theorem to vector spaces without topology, under core qualifications imposed on the doms, in line of our recent papers \cite{CBN21,CBNC}. The second branch of duality theory in optimization studies \emph{Lagrangian duality}. One of the main goals of Lagrangian duality theory is to obtain sufficient conditions which guarantee the  {\em Lagrangian strong duality}. Most of the known results in this research direction require the {\em Slater conditions}. In locally convex topological vector spaces,  Slater conditions based on generalized interiors were  introduced by Borwein and Lewis in \cite{blreport1}, and then were generalized by Jeyakumar and Wolkowicz in \cite{JW}. Further significant results were obtained in \cite{bot1,bot-revis,daniele1,daniele2} and the references therein. In this work, building on our recent papers \cite{CBN21,CBNC} on vector spaces without any topological structure, we use a Slater condition via the core of a cone to give a sufficient condition of the Lagrangian strong duality.
	
	We structure our paper in the following manner. Section~\ref{sec:preliminiaries} presents the basic notions and definitions used throughout this work. In Section~\ref{sec:generalized} some coderivative calculus rules for set-valued mappings are revisited and improved, along with a presentation of the subdifferential maximum rule. Note that the notion of coderivative is one of the most important concepts in modern convex analysis, and was introduced by B. S. Mordukhovich in the 1980s for set-valued mappings without the requirement of convexity. In Section~\ref{Fenchel_duality} we develop conjugate calculus rules and the Fenchel duality for extended-valued functions on vector spaces. Finally, Section~\ref{Lagrange_duality} provides some necessary and sufficient conditions for convex optimization including Karush-Kuhn-Tucker condition, and a Slater condition via the algebraic core  to get Lagrangian strong duality in vector spaces.
	
	\section{Preliminaries}\label{sec:preliminiaries}
	We begin by discussing the basic notions that will be used in the rest of the paper, see the books \cite{hol,ktz,mn-book,z} and the papers \cite{CBN21,CBNC} for relevant results. We use standard notation from convex and variational analysis; see, e.g., \cite{mn-book,m-book,rw,r}. Only real vector spaces are considered. Given such a vector space $X$, the algebraic dual space  of $X$ is
	\begin{equation*}
		X^\prime:=\big\{f\colon X\to\R\;\big|\ f\;\text{ is a linear function}\big\}.
	\end{equation*}
	Next, we introduce a few basic notions for partially ordered vector spaces using definitions and notations from \cite{luxemburg-1971-riesz}, more material can be found in \cite{CBNG,dtl}. Given a vector space $Y$ and a nonempty convex cone $C\subset Y$, we define a relation on $Y$ by $z_1\leq_{C}z_2\; \text{if and only if}\; z_2-z_1\in C$. This defines a {\em partial ordering} on $Y$, and we say that $Y$ is {\em partially ordered by $C$} and refer to $Y$ as an {\em ordered vector space}. Furthermore, an ordered vector space \(Y\) is considered {\em generating} if $Y=C-C$. In this paper, we only consider {\em Archimedean} ordered vector spaces, that is $n^{-1}y\downarrow_C 0$ whenever $y\in Y^+$.
	
	Now consider a function between vector spaces, $\phi\colon X\to Y$. Here, $Y$ is partially ordered by \(C\). We then say $\phi$ is {\em $C$-convex} if for $a,b\in\dom_C(\phi)$ and $\lambda\in[0,1]$
	\begin{equation*} \phi\big(\lambda a+(1-\lambda)b\big)\leq_C
		\lambda \phi(a)+(1-\lambda)\phi(b).
	\end{equation*}
	Given a set $\Lambda\subset X$, we define the {\em core} or {\em algebraic interior} of it by
	\begin{eqnarray*}\label{core}
		\co(\Lambda):=\big\{x\in\Lambda\;\big|\;\forall v\in X,\;\exists\delta>0,\;\forall t\mbox{ with }|t|<\delta:\;x+t v\in\Lambda\big\}.
	\end{eqnarray*}
	In the case $X$ is a topological vector space, we have ${\rm int}(\Lambda)\subset\co(\Lambda)\subset\Lambda$, where ${\rm int}(\Lambda)$ is the standard topological interior of $\Lambda$. A set $\Lambda\subset X$ is {\em core-solid} if $\co(\Lambda)\ne\emp$.
	
	Two nonempty subsets $\Lambda_1,\Lambda_2\subset X$ can be {\em separated} by a hyperplane if there exists a nonzero linear function $g\colon X\to\R$ such that
	\begin{equation*}\label{sep1}
		\sup\big\{g(x)\;\big|\;x\in\Lambda_1\big\}\le\inf\big\{g(x)\;\big|\;x\in\Lambda_2\big\}.
	\end{equation*}
	In addition, if $\inf\big\{g(x)\;\big|\;x\in\Lambda_1\big\}<\sup\big\{g(x)\;\big|\;x\in\Lambda_2\big\}$, that is there exists $x_1\in\Lambda_1$ and $x_2\in\Lambda_2$ with $g(x_1)<g(x_2)$, then $\Lambda_1$ and $\Lambda_2$ are {\em properly separated} by a hyperplane.
	
	In particular, the separability of a convex set $\Lambda$ and singleton $\{x_{0}\}$ can be characterized by whether $x_{0}$ exists in the core of $\Lambda$ \cite[Theorems 3.3, 3.4]{CBNC}. On the other hand, \cite[Theorem 3.2]{CBN21} and \cite[Lemma 5.2]{CBNC} gave us versions of properly separated theorems  to two convex sets in vector spaces.
	
	Let $\Lambda \subset X$ be nonempty and convex, then the {\em normal cone} to $\Lambda$ at $\oz\in\Lambda$ is given by
	\begin{equation*}\label{nc}
		N(\oz;\Lambda):=\big\{f\in X^\prime\;\big|\;f(z-\oz)\le 0\;\text{ for all }\;z\in\Lambda\big\}
	\end{equation*}
	with $N(z;\Lambda):=\emp$ for $z\notin\Lambda$. It is easy to see that if $z\in\core(\Lambda)$, then $N(z;\Lambda)=\{0\}$.  Note that \cite[Theorem 3.3]{CBN21} gave us an important rule we will use a few times,  a representation of the normal cone to the intersection of two convex sets in vector spaces.
	
	\section{Revisiting Generalized Differential Calculus in Vector Spaces}\label{sec:generalized}
	We first revisit some coderivative calculus rules with improvements utilizing qualifying conditions in \cite{CBN21} then provide the subdifferential maximum rule  in vector spaces.
	
	Consider a set-valued mapping $G\colon X\tto Y$ between vector spaces. The {\em domain}, {\em range}, and {\em graph} of it are defined by
	\begin{align*}
		& \dom(G):=\big\{x\in X\;\big|\;G(x)\ne\emp\big\},\quad \rge(G):=\bigcup_{x\in X}G(x) \\
		& \gph(G):=\big\{(x,y)\in X\times Y\;\big|\;y\in G(x)\big\}
	\end{align*}
	Recall that $G$ is convex if its graph is a convex subset of $X\times Y$, and $\rge (G)$ is also convex given the convexity of $G$.
	
	Next, we revisit the {\em coderivative}	 of $G$ at $(a,b)\in\gph(G)$ which is a set-valued mapping $D^*G(a,b)\colon Y'\tto X^\prime$ and is given by
	\begin{equation}\label{cod}
		D^*G(a,b)(h):=\big\{f\in X^\prime\;\big|\;(f,-h)\in N\big((a,b);\gph(G)\big)\big\},\quad h\in Y'.
	\end{equation}
	The definition of the {\em composition} of two set-valued mappings $G\colon X\tto Y$ and $H\colon Y\tto Z$ is given by
	\begin{equation*}
		(H\circ G)(x)=\bigcup_{y\in G(x)}H(y):=\{z\in H(y)\;|\;y\in G(x)\},\quad x\in X.
	\end{equation*}
	It is not difficult to see that $H\circ G$ is convex given the convexity of both $G$ and $H$. To derive the coderivative chain rule, let $\nu \in(H\circ G)(u)$ and we consider the set $K(u,\nu):=G(u)\cap H^{-1}(\nu)$. In the following theorem, we utilize the core qualification condition on domains and ranges, instead of the graph qualifying conditions known from \cite[Theorem 5.1]{CBN21}, to derive the coderivative chain rule for set-valued mappings in vector spaces.
	\begin{theorem}\label{scr}
		Consider convex set-valued mappings $G\colon X\tto Y$, $H\colon Y\tto Z$ and suppose that the following qualifying condition is satisfied
	\begin{equation}\label{QCCR}
	\core\big(\rge(G)\big)\cap\core\big(\dom(H)\big)\ne\emp\ \text{and}\ \core\big(\rge(G)-\dom(H)\big)\ne\emp;
\end{equation}
		Then for $(a,c)\in\gph(H\circ G)$ and $s\in Z'$ we have the coderivative chain rule
		\begin{equation}\label{chain}
			D^*(H\circ G)(a,c)(s)=\bigcap_{b\in K(a,c)}D^*G(a,b)\circ D^*H(b,c)(s).
		\end{equation}
	\end{theorem}
	{\bf Proof.} The inclusion ``$\supset$" is trivial. Let $f\in D^*(H\circ G)(a,c)(s)$ and $b\in K(a,c)$. Then, $(f,-s)\in N((a,c);\gph(H\circ G))$, which implies
	\begin{equation*}
		\la f,u-a\ra-\la s,\nu-c\ra\le 0\;\mbox{\rm for}\;(u,\nu)\in\gph(H\circ G).
	\end{equation*}
	Let us define $\Lambda_1:=\gph(G)\times Z$ and $\Lambda_2:=X\times\gph(H)$ to be convex subsets of $X\times Y\times Z$. It can be easily seen that
	\begin{equation}\label{graphrelation1}
		\Lambda_1-\Lambda_2=X\times
		\big(\rge(G)-\dom(H)\big)\times Z.
	\end{equation}
	Using \eqref{graphrelation1}, we arrive at the representation
	\begin{equation*}
		\core(\Lambda_1-\Lambda_2)= X\times\core\big(\rge(G)-\dom(H)\big)\times Z.
	\end{equation*}
	Applying \cite[Lemma 5.1]{CBNC} yields
	\begin{equation*}
		\core(\Lambda_1-\Lambda_2)= X\times\big(\core\big(\rge(G)\big)-\core\big(\dom(H)\big)\big)\times Z.
	\end{equation*}
	From the definitions of $\Lambda_1$, $\Lambda_2$ and \eqref{QCCR} we have
	\begin{equation}\label{QC1_re}
		0\in\core(\Lambda_1-\Lambda_2),\;\mbox{ and so }\;\core(\Lambda_1)\cap\core(\Lambda_2)\ne\emp.
	\end{equation}
	Observe that $(a,b,c)\in \Lambda_1\cap\Lambda_2$. Hence, from the definitions
	\begin{equation*}
		(f,0,-s)\in N((a,b,c);\Lambda_1\cap\Lambda_2).
	\end{equation*}
	Now, employing \cite[Theorem 5.4]{CBNC} with the qualifying condition \eqref{QC1_re} we obtain
	\begin{equation*}
		(f,0,-s)\in N((a,b,c);\Lambda_1\cap\Lambda_2)=N((a,b,c);\Lambda_1)+N((a,b,c);\Lambda_2),
	\end{equation*}
	and therefore we can find $k\in Y'$ such that
	\begin{equation*}
		(f,0,-s)=(f, -k,0)+(0,k,-s)
	\end{equation*}
	with $(f,-k)\in N((a,b);\gph(G)),\;(k,-s)\in N((b,c);\gph(H)).$ Using the coderivative definition (\ref{cod}) we arrive at $f\in D^*G(a,b)(k)\;\mbox{\rm and }\;k\in D^*H(b,c)(s)$, which means ``$\subset$" in \eqref{chain} is verified.$\h$
	
	Consider the given set $S\subset Y$ and a set-valued mapping $G\colon X\tto Y$. The {\em preimage} of $S$ under $G$ is denoted
	\begin{equation*}
		G^{-1}(S):=\{x\in X\ |\ G(x)\cap S\ne\emptyset\}.
	\end{equation*}
	
	The next proposition provides a formula for the normal cone to $G^{-1}(S)$ through the normal cone of $S$ and the coderivative of $G$.
	\begin{proposition}\label{Theo-code-preimage}
		Consider a convex set-valued mapping $G\colon X\tto Y$, a convex set $S\subset Y$, and assume that
\begin{equation}\label{ri_inver}
	\core\big(\rge(G)\big)\cap\core(S)\ne\emptyset\ \text{and}\ \core\big(\rge(G)-S\big)\ne\emptyset.
\end{equation}
		Then, for $u\in G^{-1}(S)$ and $v\in G(u)\cap S$ we have the formula
		\begin{equation*}\label{Con-inverse}
			N(u;G^{-1}(S))=D^*G(u,v)\big(N(v;S)\big).
		\end{equation*}
	\end{proposition}
	{\bf Proof.} Observe that the set $G^{-1}(S)$ is clearly convex. Pick $f\in N(u;G^{-1}(S))$,
	\begin{equation*}
		\la f,x-u\ra \leq 0\ \text{ for}\ x\in G^{-1}(S),\ \text{i.e.},\ G(x)\cap S\ne\emptyset.
	\end{equation*}
	Let us consider convex subsets of $X \times Y$, $ \Lambda_1:=\gph(G)\;\mbox{ and }\;\Lambda_2:=X\times S$. One can see that $\Lambda_1-\Lambda_2=X\times \big(\rge(G)-S\big).$
	Using \eqref{ri_inver}, we yield
	\begin{equation*}
		0\in \core(\Lambda_1)-\core(\Lambda_2)=\core(\Lambda_1-\Lambda_2)= X\times\big(\core(\rge(G))-\core(S)\big),
	\end{equation*}
	and hence $\core(\Lambda_1)\cap\core(\Lambda_2)\ne\emptyset$. Applying \cite[Theorem 5.4]{CBNC} we obtain
	\begin{equation*}
		\begin{array}{ll}
			(f,0)\in N((u,v);\Lambda_1\cap\Lambda_2) & =N((u,v);\Lambda_1)+N((u,v);\Lambda_2)  \\
			& = N((u,v);\gph(G))+[\{0\}\times N(v;S)]
		\end{array}
	\end{equation*}
	which gives us
	\begin{equation*}
		(f,0)=(f,-k)+(0,k)\ \text{where}\ (f,-k)\in N((u,v);\gph(G))\ \text{and}\ k\in N(v;S).
	\end{equation*}
	Hence, $f\in D^*G(u,v)(k)$ and we get the inclusion
	\begin{equation*}
		N(u;G^{-1}(S))\subset D^*G(u,v)\big(N(v;S)\big).
	\end{equation*}
	The opposite inclusion is trivial, and thus we are done.
	$\h$
	
	Next, we turn our attention into {\em parametric constraint systems} and present the precise formula for the coderivative of those. Given vector spaces $X, Y$ and $Z$, we consider convex sets $\Lambda\subset X\times Y$ and $S\subset Z$, and a mapping $\phi\colon X\times Y\to Z$ such that $\phi^{-1}(S)$ and $\rge(\phi)$ are convex. Let the set-valued mapping $G\colon X\tto Y$ be defined by
	\begin{equation}\label{G}
		G(x):=\{y\in Y\ |\ (x,y)\in \Lambda,\ \phi(x,y)\in S.\}
	\end{equation}
	Note that the graph of $G$ here is convex. We have the following result.
	
	\begin{theorem}
		Consider sets $S, \Lambda$ and the function $\phi$ as described above. Let $G$ be the set-valued mapping in \eqref{G}. Assume that $\core\big(\rge(\phi)\big)\cap \core(S)\ne\emptyset$. Furthermore, if either $\core\big(\phi^{-1}(S)\big)\cap \Lambda\ne\emptyset\; \text{or}\; \phi^{-1}(S)\cap\core(\Lambda)\ne\emptyset$ is satisfied.
		Then for $(u,v)\in \gph(G)$ with $w=\phi(u,v)\in S$ and any $h\in Y'$
		\begin{equation*}\label{code_constraint}
			D^*G(u,v)(h)=\{f\in X'\ |\ (f,-h)\in D^*\phi(u,v,w)(N(w;S))+N((u,v);\Lambda)\}.
		\end{equation*}
	\end{theorem}
	{\bf Proof.} We first notice that $\gph(G)=\phi^{-1}(S)\cap \Lambda$. Given $(u,v)\in \gph(G)$ and $h \in Y'$ we obtain the coderivative of \eqref{G} as follows
	\begin{equation*}
		D^*G(u,v)(h)=\{f\in X' \ |\ (f,-h)\in N((u,v);\phi^{-1}(S)\cap \Lambda)\}.
	\end{equation*}
	Using the assumptions and applying \cite[Theorem 3.3]{CBN21} give us
	\begin{equation*}
		N((u,v);\gph(G))=N((u,v);\phi^{-1}(S))+N((u,v);\Lambda).
	\end{equation*}
	Next, observe that $N\big((u,v);\phi^{-1}(S)\big)=D^*\phi(u,v,w)(N(w;S))$ by using Proposition \ref{Theo-code-preimage},
	and hence the proof is done.
	$\h$
	
	Given set-valued mappings $G_1,G_2\colon X\tto Y$, the intersection $G_1\cap G_2\colon X\tto Y$ is define by
	\begin{equation*}
		(G_1\cap G_2)(x):=G_1(x)\cap G_2(x),\; x\in X.
	\end{equation*}
	In the following theorem we present a representation of the coderivative for the intersection of two set-valued mapping.
	
	\begin{theorem}\label{thm:int}
		Consider convex set-valued mappings $G_1,G_2\colon X\rightrightarrows Y$ and assume
		\begin{equation*}
			\co(\gph G_1)\cap\gph G_2\ne\emp.
		\end{equation*}
		Then, for $v\in(G_1\cap G_2)(u)$ and any $g\in Y'$ we have the representation
		\begin{equation}\label{cod-inter}
			D^*(G_1\cap G_2)(u,v)(g)=\bigcup_{g_1+g_2=g}\big(D^*G_1(u,v)(g_1)+D^*G_2(u,v)(g_2)\big).
		\end{equation}
	\end{theorem}
	{\bf Proof.} The inclusion ``$\supset$'' in \eqref{cod-inter} is straightforward. We  now prove the inclusion ``$\subset$" in \eqref{cod-inter}. Let $v\in(G_1\cap G_2)(u)$, and $f\in D^*(G_1\cap G_2)(u,v)(g)$ we have
	\begin{equation*}
		(f,-g)\in N\big((u,v);{\gph}(G_1\cap G_2)\big)=N\big((u,v);(\gph G_1)\cap(\gph G_2)\big).
	\end{equation*}
	Then applying \cite[Theorem 3.3]{CBN21} under the corresponding assumptions yields
	\begin{align*}
		(f,-g)\in N\big((u,v);{\rm gph}(G_1\cap G_2)\big)=N\big((u,v);\gph G_1\big)+N\big((u,v);\gph G_2\big).
	\end{align*}
	Thus, $(f,-g)=(f_1,-g_1)+(f_2,-g_2)$ with $(f_i,-g_i)\in N((u,v);\gph G_i))$ for $i=1,2$. Therefore, $f\in D^*G_1(u,v)(g_1)+D^*G_2(u,v)(g_2)$ and $g=g_1+g_2$, which justify the inclusion ``$\subset$" in \eqref{cod-inter}. $\h$

		By a simple proof, the following lemma gives a representation for $\core(\epi(\phi))$, where $\phi\colon X\to\Bar{\R}$ is an extended real value function, which will be used in proving the next theorem.
\begin{lemma}\label{lm_epi_core_solid} Consider a convex function $\phi\colon X\to\bar{\R}$, if $\dom(\phi)$ is core-solid, then $\epi(\phi)$ is also core-solid, and we have the representation
\begin{equation*}
  \core\big(\epi(\phi)\big)=\{(x,\lambda)\in X\times \R\; \big|\; x\in \core\big(\dom(\phi)\big), \phi(x)<\lambda\}.
\end{equation*}
\end{lemma}
{\bf Proof.}
	Let $\overline{x} \in \core \left( \dom \phi \right) $. We will show that $(\overline{x},\phi(\overline{x})+1) \in \core \left( \epi (\phi) \right) $. To prove this we shall first detour and prove an auxiliary result. Define the
set
\[
A \coloneqq \{x \in X : \; \phi(x) < \phi(\overline{x}) + 1/2\}.
\]
We shall show that $A$ itself is core-solid, with $\overline{x} \in \core (A)$. Take $v \in X$. Clearly, $\overline{x} \in A$ and as $\overline{x} \in \core(\dom \phi) $ there exists $\delta_{0}$ such that for all $\abs{t }< \delta_{0}, \; \overline{x} + vt \in \dom(\phi) $. Thus the function
\[
\hat{\phi} (t)\coloneqq \phi(\overline{x}+vt), \quad t \in (-\delta_{0},\delta_{0})
\]
is real-valued and convex being the composition of a convex and affine function. We handle the case $t\geq0$ first. If $\hat{\phi}(t) < \phi(\overline{x})+ 1 / 2$ for all $0\le t< \delta_{0}$ we are already done, so suppose  not. Without loss of generality suppose that
\[
\gamma \coloneqq \max_{0\le t\le \delta_{0}} \hat{\phi}(t)
\]
is greater than $\phi(\overline{x})$. As $\hat{\phi}$ is convex, $\gamma$ either is finite, or can be made finite by reducing $\delta_{0}$, and is attained at some point $t^{*}$. We can then find the intersection of the segment connecting $\hat{\phi}(0)=\phi(\overline{x})$ and $\hat{\phi}(t^{*})$ with $\phi(\overline{x})+ 1/2$, i.e solve
\[
z \hat{\phi}(t^{*}) + (1-z)\phi(\overline{x}) = \phi(\overline{x})+1/2, \quad z \in [0,1]
\]
which yields
\[
z = \frac{1}{2 \left( \hat{\phi}(t^{*})-\phi(\overline{x}) \right) }
\]
Thus letting $\delta_{R}= zt^{*}$, for $0\le t<\delta_{R}$ by the convexity of $\hat{\phi}$
\[
\phi(\overline{x}+vt) = \hat{\phi}(t)<\hat{\phi}(\overline{x})+1/2,
\]
and consequently $\overline{x}+vt \in A$. Now we can perform a similar argument for $t<0$ to obtain a $\delta_{L }$. So for $\abs{t }< \delta_{0}=\min \{|\delta_{L }|, \delta_{R}\} $ we obtain $\overline{x}+vt \in A.$ As $v $ is arbitrary we have $\overline{x} \in \core A.$

With this result, we can now prove the lemma itself. Take $\left( v, \lambda \right) \in X \times \mathbb{R}$. Then by the above result there exists $\delta_{0}$ such that for $\abs{t }<\delta_{0}$ we have $\overline{x} +vt \in A$. Simply choosing $\delta = \min \left( \delta_{0}, \frac{1}{2\abs{\lambda}} \right)$ for $\lambda \neq 0$ we have that for $\abs{t}<\delta$
\[
\phi(\overline{x})+1+\lambda t \ge \phi(\overline{x})+1/2 > \phi(\overline{x}+vt)
\]
 with the same result trivially holding when $\lambda =0$. Therefore $\left( \overline{x}+vt, \phi(\overline{x})+1+\lambda t  \right) \in \epi (\phi)$ and thus $\left( \overline{x}, \phi(\overline{x})+1 \right) \in \core \left( \epi \phi \right) $.
Define further the set-valued mapping $G:X\tto Y$ by $G(x):=[\phi(x),\infty]$. We can easily see that $\dom(G)=\dom(\phi)$ and $\gph(G)=\epi(\phi)$ which gives us
\[\core(\gph(G))=\core(\epi(\phi))\neq\emptyset\]
Applying \cite[Theorem 4.2]{CBNC} we arrive at the conclusion of this lemma.
 $\h$

	We now derive a formula for the subdifferential of the maxima of finitely many convex functions, which will be used in proving the necessary optimality conditions for the constrained optimization problem on vector spaces, Theorem~\ref{thm: Neccessary Optimality Condition}. This will again use the normal cone intersection rule \cite[Theorem 5.4]{CBNC}. Given functions $\phi_i\colon X\to\Bar{\R}$, the maximum function $\phi\colon X\to\Bar{\R}$ is
	\begin{equation}\label{3.68}
		\phi(x):=\max\{\phi_i(x)\; \big|\; i=1,\cdots, k\},\; x\in X.
	\end{equation}
	\begin{theorem}\label{sub_max_function}
		Consider the maximal function \eqref{3.68}, where each $\phi_i\colon X\to\Bar{\R}$ is convex.  Define
		\begin{equation}\label{3.69}
			I(\bar{u}):=\{i=1,\cdots, k\; \big|\; \phi_i(\bar{u})=\phi(\bar{u})\}
		\end{equation}
		as the active index set at $\bar{u}\in \cap_{i=1}^k\dom(\phi_i)$, it always holds
		\begin{equation}\label{3.70}
			\partial \phi(\bar{u})\supset {\rm co}\big(\cup_{i\in I(\bar{u})}\partial\phi_i(\bar{u})\big).
		\end{equation}
		If in addition, the core qualification
		\begin{equation}\label{3.70_1}
			\bar{u}\in \cap_{i=1}^k\core\big(\dom(\phi_i)\big),
		\end{equation}
		is satisfied, the subdifferential maximum rule holds
		\begin{equation}\label{3.71}
			\partial \phi(\bar{u})= {\rm co}\big(\cup_{i\in I(\bar{u})}\partial\phi_i(\bar{u})\big).
		\end{equation}
	\end{theorem}
	{\bf Proof.} We first show that $\partial \phi_i(\bar{u})\subset \partial \phi(\bar{u})\; \text{holds whenever}\; i\in I(\bar{u}).$ Let $i\in I(\bar{u})$ and $f\in \partial \phi_i(\bar{u})$. Using \eqref{3.69} we obtain
	\begin{equation*}
		\la f,z-\bar{u}\ra \leq \phi_i(z)-\phi_i(\bar{u})=\phi_i(z)-\phi(\bar{u})\leq \phi(z)-\phi(\bar{u})\; \forall z \in X.
	\end{equation*}
	Hence $f\in \partial \phi(\bar{u})$ and \eqref{3.70} is deduced from the convexity of $\partial \phi(\bar{u})$.
	
	Next, we justify the inclusion $``\subset "$ in \eqref{3.71}, notice that \eqref{3.68} implies
	\begin{equation*}
		\epi(\phi)=\cap_{i=1}^k\epi (\phi_i).
	\end{equation*}
	Picking any $f\in \partial \phi(\bar{u})$, we get
	\begin{equation*}
		(f,-1)\in N\big((\bar{u},\phi(\bar{u})); \epi(\phi)\big)=N\big((\bar{u},\phi(\bar{u})); \cap_{i=1}^k\epi(\phi_i)\big).
	\end{equation*}
	Under the qualification \eqref{3.70_1}, Lemma \ref{lm_epi_core_solid} shows that $\cap_{i=1}^k\core\big(\epi (\phi_i)\big)\ne\emptyset$. Applying \cite[Theorem 5.4]{CBNC} we get
	\begin{equation*}
		(f,-1)\in N\big((\bar{u},\phi(\bar{u})); \cap_{i=1}^k\epi(\phi_i)\big) = \sum_{i=1}^k N\big((\bar{u},\bar{\gamma}); \epi(\phi_i)\big)\; \text{with}\; \bar{\gamma}:=\phi(\bar{u}).
	\end{equation*}
	Observe that $\bar{\gamma}>\phi_i(\bar{u})$ when $i\notin I(\bar{u})$, and hence $(\bar{u},\bar{\gamma})\in \core\big(\epi(\phi_i)\big)$ for such $i$, and therefore $N\big((\bar{u},\bar{\gamma}); \epi(\phi_i)\big)=\{0,0\}$ for all $i\notin I(\bar{u})$.
	Therefore,
	\begin{equation*}
		(f,-1)\in \sum_{i\in I(\bar{u})} N\big((\bar{u},\phi_i(\bar{u}); \epi(\phi_i))\big),
	\end{equation*}
	and we yield the representation
	\begin{equation*}
		(f,-1)=\sum_{i\in I(\bar{u})}(f_i,-\gamma_i),\; \text{where}\;\ (f_i,-\gamma_i)\in N\big((\bar{u},\phi(\bar{u})); \epi(\phi_i)\big).
	\end{equation*}
	We can also obtain $\gamma_i\geq 0$ and $f_i \in \gamma_i \partial\phi_i(\bar{u})$  when $i\in I(\bar{u})$.  Note that since $\sum_{i\in I(\bar{u})}\gamma_i=1$, we can rewrite $f$ as $f=\sum_{i\in I(\bar{u})}f_i\in {\rm co}\big(\cup_{i\in I(\bar{u})}\partial\phi_i(\bar{u})\big).$ $\h$
	
	\section{Conjugate Calculus and Duality in Vector Spaces}\label{Fenchel_duality}
	Here, we further develop calculus rules for Fenchel conjugates in the vector space setting, using the relaxation of qualification conditions under the convex set extremality, see more in \cite{CBN21}. These results are then used to develop a new duality theorem in vector spaces. Similar work for polyhedral settings in LCTVS can be found in \cite{CBNG22}.
	
	Let us first recall the {\em support function} definition. Consider a nonempty subset $\Lambda \subset X$, the {\em support function} that  is associated with $\Lambda$, $\sigma_{\Lambda}\colon X'\to\bar{\R}$ is defined by
	\begin{equation*}
		\sigma_{\Lambda}(f):=\sup\{\la f,u\ra\;\big|\;u\in\Lambda\},\; f\in X'.
	\end{equation*}
	In the following theorem we establish a representation for the support function of the intersection of two convex sets in general vector spaces.
	\begin{theorem}\label{sigma intersection rule}
		Consider nonempty and convex subsets $\Lambda_1, \Lambda_2 \subset X$ and assume that	$\co(\Lambda_1)\cap \Lambda_2\neq \emptyset$. Then, we have  the support function of the intersection $\Lambda_1\cap\Lambda_2$ as follows
		\begin{equation}\label{supp1}
			\big(\sigma_{\Lambda_1\cap\Lambda_2}\big)(f)=\big(\sigma_{\Lambda_1}\s\sigma_{\Lambda_2}\big)(f)\;\;\mbox{ for all }\;f\in X'.
		\end{equation}
		In addition, for any $f\in\dom(\sigma_{\Lambda_1\cap\Lambda_2})$, there exist $f_1,f_2\in X'$ such that $f=f_1+f_2$ and
		\begin{equation}\label{supp2}
			(\sigma_{\Lambda_1\cap\Lambda_2})(f)=\sigma_{\Lambda_1}(f_1)+\sigma_{\Lambda_2}(f_2).
		\end{equation}
	\end{theorem}
	{\bf Proof.} We first prove the inequality ``$\ge$" in \eqref{supp1}. Let $f\in\dom(\sigma_{\Lambda_1\cap\Lambda_2})$, and denote $\omega:=\sigma_{\Lambda_1\cap\Lambda_2}(f)\in\R$. From the definition, we have $ \la f,x\ra\le\omega\;\mbox{\rm for}\;x\in \Lambda_1\cap\Lambda_2$. We then define the nonempty convex subsets of
	$X\times\R$ by
	\begin{eqnarray*}\label{theta1}
		\begin{array}{ll}  \Pi_1:=\big\{(x,\gamma)\in
			X\times\R\;\big|\;x\in \Lambda_1,\;\gamma\le\la f,x\ra-\omega\big\}, \quad
			\Pi_2:=\Lambda_2\times[0,\infty).
		\end{array}
	\end{eqnarray*}
	From the construction of $\Pi_1$, we can easily see that
	\begin{equation*}
		\co(\Pi_1)=\{x\in \co(\Lambda_1),\;\gamma < \la f,x\ra-\omega\}
	\end{equation*}
	and hence $\co(\Pi_1)\cap \Pi_2=\emptyset$. Applying \cite[Theorem 3.2]{CBN21} we obtain the element $(0,0)\ne(g,\eta)\in X'\times\mathbb{R}$ for which \begin{equation}\label{sep f}
		\la g,x\ra+\eta\gamma_1\le\la g,y\ra+\eta\gamma_2\;\mbox{ whenever }\;(x,\gamma_1)\in \Pi_1,\;(y,\gamma_2)\in\Pi_2.
	\end{equation}
	Furthermore, there exist $(\tilde{x},\tilde{\gamma}_1)\in\Pi_1$ and $(\tilde{y},\tilde{\gamma}_2)\in\Pi_2$ satisfying
	\begin{equation}\label{sepf 1}
		\la g,\tilde{x}\ra+\tilde{\gamma}_1\eta<\la g,\tilde{y}\ra+\tilde{\gamma}_2\eta.
	\end{equation}
	If $\eta=0$ we can reduce \eqref{sepf 1} to $\la g,\tilde{x}\ra<\la g,\tilde{y}\ra$. Using $(\bar{u},\la f,  \bar{u} \ra -\omega) \in \Pi_1$ and $(\bar{u},0) \in \Pi_2$ in \eqref{sep f} gives us $\eta\geq 0$.  Now we will show that $\eta >0$. By contradiction,  suppose $\eta=0$. From \eqref{sep f} and \eqref{sepf 1}, $\Lambda_1$ and $\Lambda_2$ can be properly separated which contradicts $\co(\Lambda_1)\cap \Lambda_2\neq \emptyset$. Thus, $\eta >0$.
	
	Choosing the pairs  $(x,\langle f,x\rangle-\omega)\in \Pi_1$ in \eqref{sep f} and $(y,0)\in\Pi_2$ yields
	\begin{equation*}
		\la	g,x\ra+\eta(\la f,x\ra-\omega)\leq \la g,y\ra\;\;\mbox{\rm with }\;\eta>0,
	\end{equation*}
	which leads us to $\big\la g/\eta+f,x\big\ra+\big\la -g/\eta,y\big\ra\leq \omega$ for $x\in \Lambda_1,\;y\in\Lambda_2$. Letting $f_1:=g/\eta+f$ and $f_2:=-g/\eta$, we obtain the inequality ``$\ge$" in \eqref{supp1}.
	
	The proof of the inequality ``$\leq$'' in \eqref{supp1} is trivial and hence \eqref{supp2} follows. $\h$
	
	Recall that the {\em Fenchel conjugate} of $\phi\colon X\to\oR$ is denoted $\phi^*\colon X'\to [-\infty,\infty]$ and is given by
	\begin{equation*}
		\phi^*(f):=\sup\{\la f,u\ra-\phi(u)\ |\ u\in X\},\; f\in X'.
	\end{equation*}
	Furthermore, if $\dom(\phi)\ne \emptyset$, then $\phi^*\colon X'\to(-\infty,\infty]$ is convex and
	\begin{equation*}
		\phi^*(f):=\sup\{\la f,u\ra-\phi(u)\ |\ u\in \dom(\phi)\},\; f\in X'.
	\end{equation*}
	Now there is a connection between the Fenchel conjugate of a function and the support function of its epigraph which our geometric approach to conjugate calculus relies upon. This crucial connection is shown in the following lemma.
	
	\begin{lemma}\label{Fenchel epi} Given a proper function $\phi\colon X\to\oR$, it holds that
		\begin{equation*}
			\phi^*(f)=\sigma_{{\rm epi}(\phi)}(f,-1)\;\mbox{ whenever }\;f\in X'.
		\end{equation*}
	\end{lemma}
	Notice that the proof of this lemma is similar to ~\cite[Lemma 6.2]{CBNG22} and hence omitted here.
	
	We are now ready to present the {\em conjugate sum rule} in vector spaces.
	\begin{theorem}\label{Fenchel sum rule}
		Consider proper convex functions $\phi_1, \phi_2\colon X\to\oR$  on a vector space $X$. Suppose $\co\big(\dom(\phi_1)\big)\cap \dom(\phi_2)\neq \emptyset$. We then have the conjugate sum rule
		\begin{equation}\label{Fenchelsum} (\phi_1+\phi_2)^*(f)=\big(\phi_1^*\s \phi_2^*\big)(f)\;\mbox{ for all }\;f\in X'.
		\end{equation}
		Furthermore, the infimum in $(\phi_1^*\s \phi_2^*)(f)$ is attained, i.e., for $f\in\dom(\phi_1+\phi_2)^*$ there exist $f_1,f_2\in X'$ such that
		\begin{eqnarray*}\label{inf-conj}
			(\phi_1+\phi_2)^*(f)=\phi_1^*(f_1)+\phi_2^*(f_2),\quad f=f_1+f_2.
		\end{eqnarray*}
	\end{theorem}
	{\bf Proof.} Let $f_1,f_2\in X'$ with $f_1+f_2=f$, we get
	\begin{eqnarray*}
		\begin{array}{ll}
			\phi_1^*(f_1)+\phi_2^*(f_2) & =\sup\big\{\la f_1,z\ra-\phi_1(z)\;\big|\;z\in X\big\}+\sup\big\{\la f_2,z\ra-\phi_2(z)\;\big|\;z\in X\big\} \\
			& \ge\sup\big\{\la f_1,z\ra-\phi_1(z)+\la f_2,z\ra-\phi_2(z)\;\big|\;z\in X\big\}                              \\
			& =\sup\big\{\la f,z\ra-(\phi_1+\phi_2)(z)\;\big|\;z\in X\big\}=(\phi_1+\phi_2)^*(f).
		\end{array}
	\end{eqnarray*}
	Note that the given assumption is not used in the proof of this inequality.
	
	Next, we prove that $(\phi_1^*\s \phi_2^*)(f)\le(\phi_1+\phi_2)^*(f)$. It suffices to only consider the case when $(\phi_1+\phi_2)^*(f)<\infty$. Define two convex sets by
	\begin{eqnarray}\label{omega-conj}
		\begin{array}{ll}
			& \Lambda_1:=\big\{(z,\gamma_1,\gamma_2)\in X\times\R\times\R\;\big|\;\gamma_1\ge \phi_1(z)\big\}=\epi(\phi_1)\times\R, \\
			& \Lambda_2:=\big\{(z,\gamma_1,\gamma_2)\in X\times\R\times\R\;\big|\;\gamma_2\ge \phi_2(z)\big\}.
		\end{array}
	\end{eqnarray}
	Similarly to Lemma~\ref{Fenchel epi} we obtain the representation
	\begin{equation}\label{conj-supp}
		(\phi_1+\phi_2)^*(f)=\sigma_{\Lambda_1\cap\Lambda_2}(f,-1,-1).
	\end{equation}
	Since $\Lambda_1$ and $\Lambda_2$ satisfy the requirements in Theorem~\ref{sigma intersection rule}, applying it to \eqref{conj-supp} produces triples $(f_1,-u_1,-u_2)\in X'\times\R\times\R$ and $(f_2,-\upsilon_1,-\upsilon_2)\in X'\times\R\times\R$ such that $(f,-1,-1)=(f_1,-u_1,-u_2)+(f_2,-\upsilon_1,-\upsilon_2)$. Thus, \eqref{conj-supp} becomes
	\begin{equation*}
		(\phi_1+\phi_2)^*(f)=\sigma_{\Lambda_1\cap\Lambda_2}(f,-1,-1)=\sigma_{\Lambda_1}(f_1,-u_1,-u_2)+\sigma_{\Lambda_2}(f_2,-\upsilon_1,-\upsilon_2).
	\end{equation*}
	By the construction of $\Lambda_1$ and $\Lambda_2$, it must hold that $u_2=\upsilon_1=0$. Applying Lemma~\ref{Fenchel epi} and considering \eqref{omega-conj} gives us
	\begin{eqnarray*}
		\begin{array}{ll}
			(\phi_1+\phi_2)^*(f) & =\sigma_{\Lambda_1\cap\Lambda_2}(f,-1,-1)=\sigma_{\Lambda_1}(f_1,-1,0)+\sigma_{\Lambda_2}(f_2,0,-1) \\
			& =\sigma_{{\rm epi}(\phi_1)}(f,-1)+\sigma_{{\rm epi}(\phi_2)}(f_2,-1)                                \\
			& =\phi_1^*(f_1)+\phi_2^*(f_2)\ge\big(\phi_1^*\s \phi_2^*\big)(f).
		\end{array}
	\end{eqnarray*}
	This validates the sum rule \eqref{Fenchelsum} and hence the infimum is attained. $\h$
	
	We next introduce the conjugate chain rule, which will be used in the proof of the composite convex Fenchel duality result.
\begin{theorem}\label{Fenchel chain rule}
		Consider a linear mapping between vector spaces $X,Y$, $A\colon X\to Y$ and a convex function $\phi:Y\to \oR$. If
			\begin{equation}\label{eq1_composi_fen1}
				AX \cap \co(\dom \phi) \ne\emptyset,
			\end{equation}
		then we have a representation of the conjugate chain rule
		\begin{equation*}
			(\phi\circ A)^*(f)=\inf\{\phi^*(s)\; \big|\; s\in (A^*)^{-1}(f)\},\; f\in X'.
		\end{equation*}
		In addition, for $f\in\dom(\phi\circ A)^*$ there exists $s\in (A^*)^{-1}(f)$ such that
		\begin{equation*}
			(\phi\circ A)^*(f)=\phi^*(s).
		\end{equation*}
\end{theorem}
	{\bf Proof.} Let $s\in(A^*)^{-1}(f)$. By definition,
	\begin{eqnarray*}
		\begin{array}{ll}
			\phi^*(s) & =\sup\big\{\la s,v\ra-\phi(v)\;\big|\;v\in Y\big\}                              \\
			& \ge\sup\big\{\la s,Au\ra-\phi(Au)\;\big|\;u\in X\big\}                          \\
			& =\sup\big\{\la A^*s,u\ra-(\phi\circ A)(u)\;\big|\;u\in X\big\}                  \\
			& =\sup\big\{\la f,u\ra-(\phi\circ A)(u)\;\big|\;u\in X\big\}=(\phi\circ A)^*(f).
		\end{array}
	\end{eqnarray*}
	Therefore, $\inf\big\{\phi^*(s)\;\big|\;s\in(A^*)^{-1}(f)\big\}\ge(\phi\circ A)^*(f)$.  This inequality is also clearly true if $(A^*)^{-1}(f)=\emptyset$.
	
	Next, we prove the opposite inequality. Let $f\in\dom(\phi\circ A)^*$ and consider the following convex sets
	\begin{equation}\label{omega-chain}
		\Lambda_1:=\gph(A)\times\R \quad \text{and}\quad \Lambda_2:=X\times\epi(\phi)\subset X\times Y\times\R.
	\end{equation}
	From the structure of $\Lambda_1$ and $\Lambda_2$ one can see that
	\begin{equation*}
		(\phi\circ A)^*(f)=\sigma_{\Lambda_1\cap \Lambda_2}(f,0,-1)<\infty.
	\end{equation*}
Applying Lemma \ref{lm_epi_core_solid}, qualification \eqref{eq1_composi_fen1} gives
	\begin{equation*}
		\co(\epi \phi)=\{(y,\gamma)\; \big|\; y\in\co(\dom \phi),\;\phi(y)<\gamma\},
	\end{equation*}
	and therefore, $\Lambda_1\cap \co(\Lambda_2)\ne\emptyset$. Using Theorem~\ref{sigma intersection rule}, there exist triples $(f_1,g_1,u_1)$ and
	$(f_2,g_2,u_2)$ in $X'\times Y'\times\R$ satisfying
	\begin{equation*}
		(f,0,-1)=(f_1,g_1,u_1)+(f_2,g_2,u_2),
	\end{equation*}
	and $\sigma_{\Lambda_1\cap\Lambda_2}(f,0,-1)=\sigma_{\Lambda_1}(f_1,g_1,u_1)+\sigma_{\Lambda_2}(f_2,g_2,u_2)$. From the structure of $\Lambda_1$ and $\Lambda_2$ in \eqref{omega-chain}, it follows that $u_1=0$ and $f_2=0$. Hence, we obtain
	\begin{equation*}
		\sigma_{\Lambda_1\cap\Lambda_2}(f,0,-1)=\sigma_{\Lambda_1}(f,g_2,0)+\sigma_{\Lambda_2}(0,-g_2,-1)
	\end{equation*}
	for some $g_2\in Y'$. Thus,
	\begin{eqnarray*}
		\begin{array}{ll}
			\sigma_{\Lambda_1\cap\Lambda_2}(f,0,-1) & =\sup\big\{\la f,u\ra-\la g_2,Au\ra\;\big|\;u\in X\big\}+\sigma_{{\rm epi}(\phi)}(g_2,-1) \\
			& =\sup\big\{\la f-A^*g_2,u\ra\;\big|\;u\in X\big\}+\phi^*(g_2),
		\end{array}
	\end{eqnarray*}
	so we can conclude $f=A^*g_2$. Therefore
	\begin{equation*}
		\sigma_{\Lambda_1\cap\Lambda_2}(f,0,-1)=\phi^*(g_2)\ge\inf\big\{\phi^*(s)\;\big|\;s\in(A^*)^{-1}(f)\big\}.
	\end{equation*}
	Hence, our proof is done. $\h$
	
	The version of Fenchel strong duality in vector spaces under the core qualifications is then shown in the following theorem.
	\begin{theorem}\label{Theo_Fen_Poly} Consider convex functions $\phi,\psi \colon X\to \oR$ with further assumption $\co\left(\dom(\phi)\right) \cap \dom(\psi)\ne\emptyset$. Then
		\begin{equation*}\label{duality_Theo1}
			\inf\{\phi(x)+\psi(x)\ |\ x\in X\}=\sup\{-\phi^*(-f)-\psi^*(f)\ |\ f\in X'\}
		\end{equation*}
	\end{theorem}
	{\bf Proof.} First we prove the inequality $``\leq$". It is sufficient to only consider when $\omega:=\inf\{\phi(x)+\psi(x)\; \big|\; x\in X\}\in \R$. It's easily seen that
	\begin{equation*}
		\omega=\inf\{\phi(x)+\psi(x)\; \big|\; x\in X\}=-\sup\{\la 0,x\ra-(\phi+\psi)(x)\; \big|\; x\in X\}=-(\phi+\psi)^*(0).
	\end{equation*}
	Using the conjugate sum rule from Theorem~\ref{Fenchel sum rule}, we have a $f\in X'$ such that
	\begin{equation*}
		\omega=-(\phi+\psi)^*(0)=-\phi^*(-f)-\psi^*(f).
	\end{equation*}
	So the inequality $``\leq$" is justified while $``\geq$" is straightforward. $\h$
	
	We next examine the {\em composite convex optimization} framework of Fenchel duality in general vector spaces that was discussed in another setting in \cite{CBNG22}. Given vector spaces $X,Y$, proper convex functions $\phi\colon X\to\oR$,  $\psi\colon Y\to\oR$ and a linear  operator $A\colon X\to Y$, we can consider the following primal minimization problem:
	\begin{eqnarray}\label{FP}
		\mbox{ minimize }&\phi(x)+\psi(Ax)\;\mbox{ subject to }\;x\in X.
	\end{eqnarray}
	The {\rm Fenchel dual problem} of \eqref{FP} is
	\begin{eqnarray}\label{FD}
		\mbox{ maximize }&-\phi^*(A^*g)-\psi^*(-g)\;\mbox{ subject to }\;g\in Y'.
	\end{eqnarray}
	where $\phi^*,\psi^*$ are conjugate functions, and $A^*$ denotes the adjoint operator.
	
	We now establish a relationship between the solutions of the dual problem \eqref{FD} and the primal problem \eqref{FP}, often termed {\em weak duality}. This requires no assumption of convexity on $\phi$ or $\psi$, and the proof follows directly from definitions.
	
	\begin{proposition}\label{wd}
		Consider the primal problem \eqref{FP} and its dual problem \eqref{FD} with $\phi$ and $\psi$ not necessarily convex. Denote the optimal values for these problems by
		\begin{eqnarray*}
			\begin{array}{ll}
				\disp\Hat p:=\inf_{x\in X}\big\{\phi(x)+\psi(Ax)\big\},\quad \disp\Hat d:=\sup_{g\in Y^*}\big\{-\phi^*(A^*g)-\psi^*(-g)\big\}.
			\end{array}
		\end{eqnarray*}
		Then $\Hat p\ge\Hat d$.
	\end{proposition}
	
	In the next theorem we present a version of the Fenchel-Rockafellar theorem on vector spaces under a core qualification.
   \begin{theorem}\label{sd}
			Consider the primal problem \eqref{FP} and its dual problem \eqref{FD}. Suppose that the qualification \eqref{eq1_composi_fen1} holds. In addition to the assumptions on $\phi,\psi$, and $A$ in the problem formulation, assume either
		\begin{equation*}\label{QCD}
			\dom(\psi\circ A)\cap \co\big(\dom (\phi)\big)\ne\emptyset
		\end{equation*}
		or
		\begin{equation*}\label{QCD1}
			\co\big(\dom(\psi\circ A)\big)\cap \dom (\phi)\ne\emptyset.
		\end{equation*}
		Then $\Hat p=\Hat d$ holds. Furthermore, if $\Hat p$ is finite, the	supremum in the definition of $\Hat d$ is attained.
	\end{theorem}
	{\bf Proof.}
	Applying Proposition~\ref{wd}, we just need to prove $\Hat p\le\Hat d$. It is sufficient to only consider $\Hat p\in\R$. It is clear to see that
	\begin{eqnarray*}
		\begin{array}{ll}
			\Hat p:=\disp\inf_{x\in X}\big\{\phi(x)+\psi(Ax)\big\} & =-\disp\sup_{x\in X}\big\{\la 0,x\ra-[\phi+(\psi\circ A)](x)\big\} \\
			& =-[\phi+(\psi\circ A)]^*(0).
		\end{array}
	\end{eqnarray*}
	From the conjugate sum rule in Theorem~\ref{Fenchel sum rule}, there exists $f\in X'$ such that
	\begin{equation*}
		\Hat p=-[\phi+(\psi\circ A)]^*(0)=-\phi^*(f)-(\psi\circ A)^*(-f).
	\end{equation*}
	Applying the conjugate chain rule in Theorem~\ref{Fenchel chain rule} yields $g\in Y'$ such that
	\begin{equation*}
		A^*g=-f,\;\text{and}\;(\psi\circ A)^*(-f)=\psi^*(g).
	\end{equation*}
	Choosing $u=-g$, we obtain $\Hat p=-\phi^*(A^*u)-\psi^*(-u)\le\Hat d$, which ensures the supremum in the definition of $\Hat d$ is achieved.
	$\h$
	
	\begin{corollary}
		The first qualification condition in Theorem~\ref{sd} is satisfied if we have
		\begin{equation*}
			\dom(\psi)\cap A(\co(\dom(\phi))\neq\emptyset.
		\end{equation*}
	\end{corollary}
	{\bf Proof.} Picking a $g_0\in \dom(\psi)$ with $g_0=A(u_0)$ for some value of $u_0\in \co(\dom(\phi))$ we have that $u_0\in \dom(\psi\circ A)\cap \co(\dom(\phi))$. $\h$

%

	\section{Optimality Conditions in Convex Minimization and Lagrangian Duality in Vector Spaces}\label{Lagrange_duality}
	We first address the constrained optimization problem
	\begin{equation}\label{7.3}
		\text{minimize}\; \phi(x)\; \text{subject to}\; x\in \Pi,
	\end{equation}
	where $\Pi$ is a nonempty and convex subset of a vector space $X$, and $\phi \colon X	\to \Bar{\R}$ is convex. We call $\bar{u}\in \Pi$ an {\em optimal solution} to \eqref{7.3} if
	\begin{equation*}\label{7.4} \phi(\bar{u})\leq \phi(x)\; \text{for all}\; x\in \Pi.
	\end{equation*}
	We also consider the unconstrained version of \eqref{7.3}
	\begin{equation}\label{7.5}
		\text{minimize}\; \theta(x) \; \text{on}\; X,\; \text{where}\; \theta(x):=\phi(x)+\delta_{\Pi}(x).
	\end{equation}
	Here, $\delta_{\Pi}(\cdot)$ is the indicator function associated with $\Pi$. The following lemma presenting the relationship between \eqref{7.3} and \eqref{7.5} is easy to show.
	\begin{lemma}\label{7.14}
		Given $\bar{u} \in \Pi$, in the convex setting of \eqref{7.3} the following are equivalent
		\begin{enumerate}
			\item $\bar{u}$ is an optimal solution to the constrained problem \eqref{7.3}.
			\item $\bar{u}$ is an optimal solution to the unconstrained problem \eqref{7.5}.
		\end{enumerate}
	\end{lemma}
	Now we develop basic necessary and sufficient optimality conditions for optimal solutions to \eqref{7.3} in the following theorem.
	\begin{theorem}\label{Theo7.15}
		Consider the setting of problem \eqref{7.3} with $\bar{u}\in\Pi\cap\dom(\phi)$. Suppose either
		\begin{equation}\label{7.6}
			\Pi\cap \core\big(\dom(\phi)\big)\ne\emptyset \quad \text{ or } \quad \core(\Pi)\cap\dom(\phi)\ne\emptyset.
		\end{equation}
		Then the following statements are equivalent
		\begin{enumerate}
			\item $\bar{u}$ is an optimal solution of \eqref{7.3}.
			\item $0\in\partial \phi(\bar{u})+N(\bar{u};\Pi)$.
		\end{enumerate}
	\end{theorem}
	{\bf Proof.} Applying Lemma \ref{7.14}, we only need to consider the unconstrained version \eqref{7.5} of problem \eqref{7.3}. Using \cite[Proposition 3.29]{mn-book} (the result also holds when $X$ is a vector space without topological structure)  provides us the equivalent characterization that $0\in\partial \theta(\bar{u})$. Applying the subdifferential sum rule from \cite[Corollary 4.3]{CBN21}  under the qualification conditions \eqref{7.6} yields
	\begin{equation*}
		\partial\theta(\bar{u})=\partial \phi(\bar{u})+\partial \delta_{\Pi}(\bar{u})=\partial \phi(\bar{u})+N(\bar{u};\Pi).
	\end{equation*}
	Thus, the proof is done. $\h$
	
	Next, we add a functional constraint by a $C$-convex function to a partially ordered vector space by a convex cone. Let $\psi\colon X\to Y$, where  $Y$ is an Archimedean and generating partially ordered vector space by a cone $C\subset Y$. Lastly, let $\Lambda$ be a convex subset of $X$ and $\phi\colon X\to \Bar{\R}$ a convex function. Then we will consider the problem
	\begin{equation}\tag{P} \label{PP2_2}
		\begin{aligned}
			& \text{minimize } \phi(x),                       \\
			& \text{subject to }x\in \Lambda,\;\psi(x)\in -C,
		\end{aligned}
	\end{equation}
	which can be rewritten in the form of \eqref{7.3} with
	\begin{equation}\label{PI}
		\Pi=\{x\in\Lambda\;|\; \psi(x)\in -C\}.
	\end{equation}
	Each $x\in \Pi$ is called a {\em feasible solution} of \eqref{PP2_2}. The Lagrange function for \eqref{PP2_2} is then defined by
	\begin{equation}\label{Lagrange_Function}
		L(x,f):=\phi(x)+\la f,\psi(x)\ra,\; x\in X,\; f\in Y^{'}.
	\end{equation}
	For the next two theorems we will focus on the special case when $Y=\R^m$ and $C=\R^m_+:=\{(x_1,\dots,x_m)\; \big|\; x_i\geq 0\}$. We are inspired by the work for \(X=\R^n\) by \cite{mn-easy}, and generalize the results to arbitrary vectors spaces. Then,  $\psi(x)=\big(\psi_1(x),\dots,\psi_m(x)\big)$ and for every feasible solution of \eqref{PP2_2} we can define the collection of active indices, given by
	\begin{equation*}\label{Feasible_Solution}
		I(x):=\big\{i\in\{1,\dots, m\}\; \big|\; \psi_i(x)=0\big\}.
	\end{equation*}
	First, we present necessary optimality conditions for \eqref{PP2_2} in this special case.
	\begin{theorem}\label{thm: Neccessary Optimality Condition}
		Suppose that  $Y=\R^m$ and $C=\R^m_+$. Let $\bar{u}$ be an optimal solution to \eqref{PP2_2} such that
		\begin{equation}\label{Quali_Solution_core}
			\bar{u}\in \core\big(\dom(\phi)\big)\cap\big(\cap_{i=1}^m\core(\dom(\psi_i))\big)\cap \Lambda.
		\end{equation}
		Then there exist non negative multipliers $\gamma_0,\gamma_1,\dots,\gamma_m$, not all zero such that
		\begin{equation}\label{eq_kkt1}
			0\in \gamma_0\partial \phi(\bar{u})+\sum_{i=1}^{m}\gamma_i\partial\psi_i(\bar{u})+N(\bar{u};\Lambda),
		\end{equation}
		and $\gamma_i\psi_i(\bar{u})=0$ for $i=1,\dots,m$.
	\end{theorem}
	{\bf Proof.} Let us first define the convex maximum function
	\begin{equation*}
		\eta(x):=\max\{\phi(x)-\phi(\bar{u}),\;\psi_i(x)\; \big|\; i=1,\dots,m\},\; x\in X.
	\end{equation*}
	Notice $\bar{u}$ is a solution to
	\begin{equation*}
		\text{minimize}\; \eta(x)\; \text{subject to}\; x\in \Lambda.
	\end{equation*}
	Under the qualification \eqref{Quali_Solution_core}, Theorem \ref{Theo7.15} and Theorem \ref{sub_max_function}  show that
	\begin{equation*}
		0\in {\rm co}\left(\partial \phi(\bar{u})\cup\left(\cup_{i\in I(\bar{u})}\partial \psi_i(\bar{u})\right)\right)+N(\bar{u};\Lambda).
	\end{equation*}
	This yields $\gamma_0\geq 0$, $\gamma_i\geq0$ for $i\in I(\bar{u})$ with $\gamma_0+\sum_{i\in I(\bar{u})}\gamma_i=1$ and
	\begin{equation*}
		0\in \gamma_0\partial \phi(\bar{u})+\sum_{i\in I(\bar{u})}\gamma_i\partial\psi_i(\bar{u})+N(\bar{u};\Lambda).
	\end{equation*}
	When $i\notin I(\bar{u})$, we can pick $\gamma_i=0$ to get $\gamma_i\psi_i(\bar{u})=0$, while still maintaining $(\gamma_0,\gamma_1,\dots,\gamma_m)\ne0$. Then the theorem is proven.
	$\h$
	
	The Slater condition/constraint for \eqref{PP2_2} states that
	\begin{eqnarray}\label{Slater2_2}
		\begin{array}{ll}
			\exists \Hat{x}\in \Lambda\; \mbox{\rm such that}\;
			\psi(\Hat{x})\in -\core(C).
		\end{array}
	\end{eqnarray}
	If $Y=\R^m$ and $C=\R^m_+$, then $\psi(\Hat{x})\in -\core(C)$ reduces to $\psi_i(\Hat{x})<0, \;  i=1,\dots,m$. With this we develop necessary and sufficient conditions for \eqref{PP2_2} in Karush-Kuhn-Tucker (KKT) form.
	\begin{theorem}\label{Theo_KKT}
		Consider $Y=\R^m$ and $C=\R_+^m$. Consider a feasible solution $\bar{u}$ of \eqref{PP2_2} such that
		\begin{equation*}\label{Quali_Solution_core_2}
			\bar{u}\in \core\big(\dom(\phi)\big)\cap\big(\cap_{i=1}^m\core(\dom(\psi_i))\big).
		\end{equation*}
		Suppose the Slater condition \eqref{Slater2_2} is satisfied for \eqref{PP2_2}. Then $\bar{u}$ is an optimal solution to \eqref{PP2_2} if and only if there are nonnegative Lagrange multipliers $\gamma_1,\dots,\gamma_m$ such that
		\begin{equation}\label{eq_kkt2}
			0\in \partial \phi(\bar{u})+\sum_{i=1}^{m}\gamma_i\partial\psi_i(\bar{u})+N(\bar{u};\Lambda),
		\end{equation}
		and $\gamma_i\psi_i(\bar{u})=0$ for $i=1,\dots,m$.
	\end{theorem}
	{\bf Proof.} We first prove the necessity of \eqref{eq_kkt2},  suppose by contradiction that $\gamma_0=0$ in \eqref{eq_kkt1}. Then we can find $0\ne(\gamma_1,\dots,\gamma_m)\in\R_+^m$, $f_i\in\partial\psi_i(\bar{u})$, and $f\in N(\bar{u};\Lambda)$ satisfying
	\begin{equation*}
		0=\sum_{i=1}^{m}\gamma_i f_i+f.
	\end{equation*}
	Therefore,
	\begin{equation*}
		0=\sum_{i=1}^{m}\gamma_i\la f_i,x-\bar{u}\ra +\la f,x-\bar{u}\ra\leq \sum_{i=1}^{m}\gamma_i\big(\psi_i(x)-\psi_i(\bar{u})\big)=\sum_{i=1}^{m}\gamma_i\psi_i(x)\; \text{for}\; x\in\Lambda,
	\end{equation*}
	which contradicts the Slater condition \eqref{Slater2_2}.
	
	The sufficiency of \eqref{eq_kkt2} is shown by choosing a $f_0\in\partial\phi(\bar{u}), \; f_i\in\partial\psi_i(\bar{u})$, and $f\in N(\bar{u};\Lambda)$ such that
	\begin{equation*}
		0=f_0+\sum_{i=1}^{m}\gamma_i f_i+f
	\end{equation*}
	with $\gamma_i\psi_i(\bar{u})=0$ and $\gamma_i\geq 0$ for $i=1,\dots,m$. From the definition of subdifferential and normal cone we directly obtain that $\bar{u}$ is an optimal solution of \eqref{PP2_2}. Explicitly, for $x\in\Lambda$ with $\psi_i(x)\leq 0, \; i=1,\dots,m$
	\begin{equation*}
		\begin{aligned}0 & =\la f_0+\sum_{i=1}^{m}\gamma_if_i+f,x-\bar{u}\ra =\la f_0,x-\bar{u}\ra +\sum_{i=1}^{m}\gamma_i\la f_i,x-\bar{u}\ra +\la f,x-\bar{u}\ra    \\
			& \leq \phi(x)-\phi(\bar{u})+\sum_{i=1}^{m}\gamma_i\big(\psi_i(x)-\psi_i(\bar{u})\big)=\phi(x)-\phi(\bar{u})+\sum_{i=1}^{m}\gamma_i\psi_i(x) \\
			& \leq \phi(x)-\phi(\bar{u}),
		\end{aligned}
	\end{equation*}
	which completes the proof.
	$\h$
	
	Consider the Lagrange function for \eqref{PP2_2} given by \eqref{Lagrange_Function}, the {\em Lagrange dual function} $L'\colon Y^{'}\to [-\infty,\infty)$ is defined by
	\begin{equation*}\label{Lar_dual_Func}
		L'(f)=\inf_{x\in \Lambda}\{ \phi(x)+\la f,\psi(x)\ra\},
	\end{equation*}
	and the Lagrange dual problem of \eqref{PP2_2} is given by
	\begin{eqnarray}\label{DP2_2}
		\begin{array}{ll}
			& \mbox{\rm maximize }L'(f),\; \mbox{\rm subject to }f\in C^{'},
		\end{array}
	\end{eqnarray}
	where $ C':=\{f\in Y'\; |\; 0\leq \la f,c\ra\; \text{for all}\; c\in C\}$ is the dual cone of $C$. Then the optimal values in \eqref{PP2_2} and \eqref{DP2_2} are given by
	\begin{equation*}
		p:=\inf\{\phi(x)\;|\; x\in \Pi\}, \quad
		d:=\sup\{L'(f)\;|\; f\in C'\},
	\end{equation*}
	where $\Pi$ is the set of feasible solutions of \eqref{PP2_2} defined in \eqref{PI}.
	
	We now return to the case of a general $Y$ as written in \eqref{PP2_2}.  At the end of this section we provide a sufficient condition which guarantees the  {\em  Lagrangian strong duality} ($p=d$) on  vector spaces without topological structure. Results of this kind are one of the main goals of Lagrangian duality theory.
	
	To proceed, we consider a system  of inequalities with the unknown $x\in X$:
	\begin{equation}\label{IP2_2} \left\{
		\begin{aligned}
			& \phi(x)<0,     \\
			& \psi(x)\in -C, \\
			& x\in \Lambda.
		\end{aligned}\right.
	\end{equation}
	We also consider a system of inequalities  with unknowns $h\in Y'$:
	\begin{equation}\label{ID2_2} \left\{
		\begin{aligned}
			& \phi(x)+\la h, \psi(x)\ra\geq 0 \; \mbox{\rm for all }x\in \Lambda, \\
			& h \in C'.
		\end{aligned}\right.
	\end{equation}
	These two systems are deeply connected, as shown in the next lemma.
	
	\begin{lemma}\label{lm_solu_infinite2}
		Suppose the Slater condition \eqref{Slater2_2} holds for \eqref{PP2_2}. Then \eqref{IP2_2} has a solution if and only if \eqref{ID2_2} has no solution.
	\end{lemma}
	{\bf Proof.} Assume \eqref{IP2_2} has a solution $x_0$ and suppose by contradiction \eqref{ID2_2} has a solution $\kappa_0\in C'$. Then for the solution $x_0$ of \eqref{IP2_2} one has
	\begin{equation*}
		\phi(x_0)+\la \kappa_0,\psi(x_0)\ra\leq \phi(x_0)<0,
	\end{equation*}
	a clear contradiction. Note the Slater condition \eqref{Slater2_2} is not used here.
	
	For the opposite direction, suppose \eqref{ID2_2} is not solvable and by contradiction suppose \eqref{IP2_2} isn't either. Consider the set
	\begin{equation*}
		\Theta:=\{(\eta_0, \eta)\in \R\times Y\; |\; \phi(x)< \eta_0,\; \psi(x)\leq_C\eta, \; \mbox{\rm for some }x\in \Lambda\}.
	\end{equation*}
	We first show that $\Theta$ is a nonempty like-convex set and $(0,0)\notin \Theta$. Since \eqref{ID2_2} has no solution, we can find $x_0\in\Lambda$ such that for $\kappa\in C'$ we get the strict inequality
	\begin{equation*}
		\phi(x_0)+\la \kappa,\psi(x_0)\ra <0.
	\end{equation*}
	Fix $\Bar{\kappa}\in C'$, denote $\Bar{\eta}=\psi(x_0)$ and $\Bar{\eta}_0=-\la \Bar{\kappa},\psi(x_0)\ra$. Then $(\Bar{\eta}_0,\Bar{\eta})\in \Theta$ and hence $\Theta$ is nonempty. Since \eqref{IP2_2} is not solvable, $(0,0)\notin \Theta$. In addition, the like-convexity of $\Theta$ is implied by the $C$-convexity of $\psi$ along with the convexity of $\phi$ .
	
	Next we show that $\core(\Theta)\ne \emptyset$. By the Slater condition \eqref{Slater2_2} there exists \( \hat{x} \in \Lambda\) such that \( \psi(\hat{x}) \leq_C 0\)
	and observe that
	\begin{align*}
		\Theta \supset \{ (\eta_0,\eta) \in \R \times Y\; | \; & \phi(\hat{x}) < \eta_0 , \; \psi(\hat{x}) \leq_C \eta \}                                                      \\
		& =\{\eta_0 \in \R \; | \; \phi(\hat{x}) < \eta_0 \} \times \{ \eta \in Y \; | \; \psi(\hat{x}) \leq_C \eta \}.
	\end{align*}
	Note that as \( \psi(\hat{x}) \leq_C 0\), we have \(  C \subset \{ \eta \in Y \; | \; \psi(\hat{x}) \leq_C \eta \}\). Hence,
	\begin{equation*}
		\Theta \supset \{ \eta_0 \in \R \; | \; \phi(\hat{x}) < \eta_0 \} \times C.
	\end{equation*}
	Now \( \core(\{ \eta_0 \in \R \; | \; \phi(\hat{x}) < \eta_0 \} )\) is trivially non-empty due to being an open interval in \(\R\). Therefore, using the Slater condition
	\begin{align*}
		\core(\{ \eta_0 \in \R \; | \; \phi(\hat{x}) < \eta_0 \} \times C) & = \core(\{ \eta_0 \in \R \; | \; \phi(\hat{x}) < \eta_0 \} )\times \core(C) \neq \emptyset.
	\end{align*}
	Hence \( \emptyset \neq	\core(\{ \eta_0 \in \R \; | \; \phi(\hat{x}) < \eta_0 \} \times C) \subset \core(\Theta) \).
	
	Since $(0,0)\notin \Theta$ and $\core(\Theta)\ne \emptyset$, \cite[Theorem 2.4]{CBN21} shows that $\Theta$ and $\{(0,0)\}$ can be properly separated.  This gives us multipliers  $\Tilde{\lambda}_0\in \R, \;\Tilde{\kappa}^*_0\in Y'$, not all zero, such that
	\begin{equation}\label{dual separation2_2}
		\Tilde{\lambda}_0\eta_0+\la \Tilde{\kappa}_0^*,\eta\ra\geq 0
	\end{equation}
	whenever $(\eta_0, \eta)\in \Theta$.  In addition, there exists $(\Tilde{\eta}_0,\Tilde{\eta})\in\Theta$ satisfying the strict inequality
	\begin{equation}\label{dual_sepa2_3}
		\Tilde{\lambda}_0\Tilde{\eta}_0+\la \Tilde{\kappa}^*_0,\Tilde{\eta}\ra> 0.
	\end{equation}
	By passing to the limit, observe that the inequality \eqref{dual separation2_2} also holds for any $(\eta_0,\eta)\in \Bar{\Theta}$. We will show that $\Tilde{\eta}_0> 0$ and $\Tilde{\kappa}_0^*\in C'$.
	By  the Slater condition \eqref{Slater2_2}, there exists $\Hat{x} \in\Lambda$ such that $\psi(\Hat{x})\leq_C0$.   Fix any $\eta_0>0$ such that $\eta_0>\phi(\Hat{x}),$ one has $(\eta_0,0)\in \Theta$ and \eqref{dual separation2_2} shows that $\Tilde{\lambda}_0\eta_0\geq 0.$ If $\Tilde{\lambda}_0<0,$ then $\eta_0\leq 0,$ and hence $\Hat{x}$ is a solution of \eqref{IP2_2}. This contradiction shows that $\Tilde{\lambda}_0\geq 0.$  To show $\Tilde{\kappa}_0^*\in C'$, by a contradiction, suppose that $\Tilde{\kappa}_0^*\notin C'$. Thus there exists $u_0\in C$, and hence $0\leq_Cu_0$, such that $\la \Tilde{\kappa}_0^*, u_0\ra <0$. If $\Tilde{\lambda}_0=0$, then  we can see that $(\eta_0,u_0)\in \Theta$ while \eqref{dual separation2_2} fails. If $\Tilde{\lambda}_0>0,$  set $u_1=(1-\frac{\Tilde{\lambda}_0\eta_0}{\la \Tilde{\kappa}_0^*, u_0\ra})u_0$. Since $\Tilde{\lambda}_0\eta_0\geq0$ and $\la \Tilde{\kappa}_0^*, u_0\ra <0$,  one has $u_1\in C$ and hence $\psi(\Hat{x})\leq_C0\leq_C u_1$. We can see that $(\eta_0,u_1)\in\Theta$ but  \eqref{dual separation2_2} is not satisfied. This contradiction verifies that $\Tilde{\lambda}_0\geq0$ and $\Tilde{\kappa}_0^*\in C'$.
	
	Since $\Tilde{\kappa}_0^*\in C'$ and $\psi(\Hat{x})\in -C$, one has
	\begin{equation}\label{pro-sepa-C2_2-2}
		\la \Tilde{\kappa}_0^*, \psi(\Hat{x})\ra \leq0\leq \la \Tilde{\kappa}_0^*,c\ra \; \text{for all}\; c\in C.
	\end{equation}
	Fix $x\in\Lambda$ and see that
	$(\phi(x),\psi(x))\in \Bar{\Theta}$. Thus,
	\begin{equation}\label{dual-se2_4}
		\Tilde{\lambda}_0\phi(x)+\la\Tilde{\kappa}_0^*,\psi(x)\ra \geq 0.
	\end{equation}
	We are going to show that $\Tilde{\lambda}_0>0$. By a contradiction, suppose that  $\Tilde{\lambda}_0=0$. Then, \eqref{dual-se2_4} shows $\la \Tilde{\kappa}_0^*, \psi(\Hat{x})\ra \geq 0$ and hence $\la \Tilde{\kappa}_0^*, \psi(\Hat{x})\ra = \la \Tilde{\kappa}_0^*, -\psi(\Hat{x})\ra= 0$ follows from \eqref{pro-sepa-C2_2-2}. Therefore, \eqref{pro-sepa-C2_2-2} can be written as
	\begin{equation}\label{pro-sepa-C2_2-4}
		\la \Tilde{\kappa}_0^*, -\psi(\Hat{x})\ra =0\leq \la \Tilde{\kappa}_0^*,c\ra \; \text{for all}\; c\in C.
	\end{equation}
	Using \eqref{dual_sepa2_3} gives us that $\la \Tilde{\kappa}_0^*, \Tilde{\eta}\ra >0$.
	Since $Y$ is generating partially ordered by $C$ and $\Tilde{\eta}\in Y$, {\cite[Proposition 1.2]{wong-Ng} shows that there exists $\Tilde{c}\in C$ such that $\Tilde{\eta}\leq_C\Tilde{c}$. So, $\Tilde{c}-\Tilde{\eta}\in C$  and hence
		\begin{equation}\label{pro-sepa-C2_2-1}
			\la \Tilde{\kappa}_0^*, -\psi(\Hat{x})\ra = 0<\la \Tilde{\kappa}_0^*, \Tilde{\eta}\ra\leq \la \Tilde{\kappa}_0^*, \Tilde{c}\ra.
		\end{equation}
		The inequalities \eqref{pro-sepa-C2_2-4} and \eqref{pro-sepa-C2_2-1} show that the sets $C$ and $\{-\psi(\Hat{x})\}$ can be properly separated. Hence, $\psi(\Hat{x})\notin-\core(C)$, which contradicts with the Slater condition \eqref{Slater2_2}. Therefore, $\Tilde{\lambda}_0>0$. Dividing \eqref{dual-se2_4} by $\Tilde{\lambda}_0$ shows \eqref{ID2_2} has a solution. This contradicts our assumption, and so our proof is complete.   $\h$
		
		\begin{theorem}\label{Theo-Week-Lag}
			Let $X,Y$ be vector spaces.  Let $x$ be a feasible solution of the primal problem \eqref{PP2_2} and $f$ a feasible solution of its dual problem \eqref{DP2_2}. Then  $\phi(x)\geq L'(f).$
		\end{theorem}
		{\bf Proof.} The conclusion is implied directly from the problems' definitions.
		$\h$
		
		Applying Lemma \ref{lm_solu_infinite2} and Theorem \ref{Theo-Week-Lag}, we next state and prove the Lagrangian strong duality in vector spaces.
		\begin{theorem}
			Consider the primal problem \eqref{PP2_2} and its dual problem \eqref{DP2_2}. Suppose the Slater condition \eqref{Slater2_2} is satisfied. Then the Lagrangian duality holds, i.e., $p=d$.
		\end{theorem}
		{\bf Proof.}  If $p=-\infty$, then $d=-\infty$ by Theorem \ref{Theo-Week-Lag}. Thus, we can assume that $p\in \R$. Then the system \eqref{IP2_2} can be extended to
		\begin{equation*}
			\phi(x)-p<0, \psi(x) \in -C\;, x\in \Lambda
		\end{equation*}
		and it has no solution. Thus, applying Lemma \ref{lm_solu_infinite2}, the corresponding extension system of  \eqref{ID2_2} has a solution, i.e., there exist $h\in C'$  such that
		\begin{equation*}
			\phi(x)-p+\la h, \psi (x)\ra\geq 0\; \mbox{\rm for all }x\in \Lambda.
		\end{equation*}
		This gives us
		\begin{equation*}
			d=\inf_{x\in \Lambda}\left\{\phi(x)+\la h, \psi (x)\ra\right\}\geq p.
		\end{equation*}
		Combining this inequality with Theorem \ref{Theo-Week-Lag} yields the result. $\h$
		\section{Conclusion}
		This paper showed our efforts in revisiting and improving coderivative calculus rules in general vector spaces while comprehensively utilizing  the qualifying conditions in \cite{CBN21}. We derived the subdifferential maximum rule by building on the normal cone intersection rule of \cite[Theorem 5.4]{CBNC}. In addition, we successfully developed conjugate calculus rules and the Fenchel duality for extended-valued functions on vector spaces in Section~\ref{Fenchel_duality}. Finally, Section~\ref{Lagrange_duality} provided some necessary and sufficient conditions for convex optimization including the Karush-Kuhn-Tucker
		condition in vector spaces and, under a core-solid Slater condition, the Lagrangian strong duality in vector spaces without topological structure.

	\end{document}